\documentclass[11pt,A4]{amsart}
\usepackage{ifthen,latexsym,amssymb,amsmath}
\newtheorem*{theorem}{Theorem}

\begin{document}

\author{Bart\l{}omiej Bzd\c{e}ga}

\address{Str\'o\.zy\'nskiego 15A/20 \\ 60-688 Pozna\'n, Poland}

\email{exul@wp.pl}



\title{Sparse binary cyclotomic polynomials}

\maketitle

\begin{abstract}
We derive a lower and an upper bound for the number of binary cyclotomic polynomials $\Phi_m$ with at most $m^{1/2+\varepsilon}$ nonzero terms.
\end{abstract}

\section*{Introduction}

A cyclotomic polynomial $\Phi_m \in \mathbb{Z}[x]$ is the monic polynomial of minimal degree having all the primitive $m$th roots of unity as its zeros. We say that the number $m$ and the polynomial $\Phi_m$ are \emph{binary} if $m$ is a product of two distinct odd primes.

A. Migotti \cite{Migotti-BinaryBounds} proved that a binary cyclotomic polynomial $\Phi_m$ has coefficients in $\{-1,0,1\}$ only. The explicit number $\theta_m$ of nonzero terms of $\Phi_m$ was derived by L. Carlitz \cite{Carlitz-BinaryTerms}. He proved that for $m=pq$ we have $\theta_m = 2p'q'-1$, where $q'$ denotes the inverse of $q$ modulo $p$ and similarly $p'$ is the inverse of $p$ modulo $q$.

It can be easily proved that for binary $m$ we have $m^{1/2} < \theta_m < m/2$. H.W. Lenstra proved in \cite{Lenstra-Vanishing} that for every $\varepsilon>0$ there exist infinitely many binary numbers $m$ such that $\theta_m < m^{8/13+\varepsilon}$. His method is based on the result of C. Hooley \cite{Hooley} that for every integer $a\neq0$ and every $\varepsilon>0$ there exist infinitely many primes $p$ for which $P(p-1)>p^{5/8-\varepsilon}$, where $P(n)$ denotes the largest prime factor of $n$. The constant $8/13$ is presently the best possible.

The result of C. Hooley has been improved by several authors. The best result to the date is due to R.C. Baker and G. Harman \cite{BakerHarman}, who proved that $P(p-1)>p^{0.677}$ for infinitely many primes $p$. It gives nearly $0.6$ instead of $8/13$. We cannot improve the result of R.C. Baker and G. Harman, so we present a different method to achieve our goals.

We consider the set $A_{\varepsilon}(N)$ of integers $n<N$ for which $P(n)>n^{1-\varepsilon}$ and $P(n+1)>(n+1)^{1-\varepsilon}$. By the result of A. Hildebrand \cite{Hildebrand} the set $A_{\varepsilon} = A_{\varepsilon}(\infty)$ has a positive lower density for every $\varepsilon>0$. We use this fact to prove the following theorem.

\begin{theorem}
Let $B_{\varepsilon}(N)$ denote the set of binary $m<N$ for which $\theta_m<m^{1/2+\varepsilon}$. Then we have
$$B_{\varepsilon}(N) = \left\{ \begin{array}{ll}
\Omega(N^{1/2}) & \text{for } 0 < \varepsilon < 1/2, \\
O(N^{1/2+\varepsilon}) & \text{for } 0 < \varepsilon < 1/6, \\
O(N/\log^2 N) & \text{for } 0 < \varepsilon < 1/2,
\end{array}\right.$$
where we used the $O$ and $\Omega$ asymptotical notation.
\end{theorem}

It is a well known result of Landau that
$$\#\{m\le N: m \text{ binary}\} \sim N\log\log N/\log N.$$
From this and from the third inequality it follows that for every $\varepsilon\in(0,1/2)$ the set $B_{\varepsilon} = B_{\varepsilon}(\infty)$ has relative density $0$ in the set of binary $m$.

\section*{Proof of the Theorem}

\textbf{Part I.} For every $n \in A_{\varepsilon}$ put $p=P(n)$, $q=P(n+1)$ and $m=pq$. Then
$$q'=\min\{a>0 \; : \; p\mid aq-1\} \le (n+1)/q.$$
By the definition of $A_{\varepsilon}$ we have
$$\theta_m < 2qq' \le 2n+2 < 2m^{1/(2-2\varepsilon)}+2 = 2m^{1/2 + \varepsilon/(2-2\varepsilon)}+2 < m^{1/2+\varepsilon}$$
for $m>m_0$, where $m_0$ depends only on $\varepsilon$. Moreover
$$n < m^{1/2+\varepsilon} < m,$$
hence $n=pp'-1$ or $n=qq'-1$. If both of them are in $A_{\varepsilon}$, then
$$2m^{1/2+\varepsilon} > (pp'-1)+(qq'-1) = m-1.$$
So for $m>m_1$ (where $m_1$ also depends only on $\varepsilon$) we can determine $n$ uniquely by $m$. We have $m>n$, thus every $n>m_1$ is determined uniquely by $m$.

Let $M=\max\{m_0,m_1\}$. By the inequality
$$m \le n(n+1)/2 < n^2$$
the function
$$f: A_{\varepsilon}(N^{1/2})\backslash [1,M] \to B_{\varepsilon}(N), \quad f(n)=P(n)P(n+1)$$
is injective and so
$$\# B_{\varepsilon}(N) = \Omega(\# A_{\varepsilon}(N^{1/2})) = \Omega(N^{1/2})$$
by the result of A. Hildebrand mentioned in the introduction.

\medskip\textbf{Part II.} For $m=pq\in B_{\varepsilon}$ put $n = \min\{pp',qq'\}-1$. The following facts
$$pp'+qq' = m+1, \quad (pp')(qq') < (m+1)m^{1/2+\varepsilon}$$
imply that $n < C m^{1/2+\varepsilon}$ for some constant $C$ depending only on $\varepsilon$. Also
$$m^{1/2+\varepsilon} > \theta_m=2p'q'-1 > p,q,$$
thus $p,q>m^{1/2-\varepsilon}$. We have $p \mid n$ and $q \mid n+1$ (or inversely). Moreover
$$\frac{\log n}{\log\min\{p,q\}} \le \frac{\log C+(1/2+\varepsilon)\log m}{(1/2-\varepsilon)\log m} = \frac{1/2+\varepsilon}{1/2-\varepsilon}+o(1).$$
Thus $p=P(n)$, $q=P(n+1)$ (or inversely) for $m$ large enough, because $(1/2+\varepsilon)/(1/2-\varepsilon)<2$ for $\varepsilon<1/6$. We define $m_2$ to be the smallest number satisfying the following condition:
$$\text{if } m>m_2 \text{ then } m=P(n)P(n+1).$$
It is obvious that $m_2$ depends only on $\varepsilon$. Thus the function
$$g : B_{\varepsilon}(N)\backslash[1,m_2] \to A_{\varepsilon}(CN^{1/2+\varepsilon}), \quad g(pq)=\min\{pp',qq'\}-1$$
is an injection. Finally
$$\# B_{\varepsilon}(N) = O(\# A_{\varepsilon}(CN^{1/2+\varepsilon})) = O(N^{1/2+\varepsilon})$$
due to the result of A. Hildebrand.

\medskip\textbf{Part III.} We assume that $m=pq$ with $q>p$. By the inequality $\theta_m>q$, if $\theta_m<m^{1/2+\varepsilon}$ then $q<p^{(1/2+\varepsilon)/(1/2-\varepsilon)}$. Thus
\begin{eqnarray*}
B_{\varepsilon}(N) & = & O \left( \sum_{p<N^{1/2-\varepsilon}}\left(\pi(p^{(1/2+\varepsilon)/(1/2-\varepsilon)})-\pi(p)\right)\right. \\
& & \left. + \sum_{N^{1/2-\varepsilon}<p<N^{1/2}}\left(\pi(N^{1/2})-\pi(p)\right)\right) \\
& = & O \left( \pi(N^{1/2-\varepsilon})\pi(N^{1/2+\varepsilon}) + \pi(N^{1/2})^2 \right) = O(N/\log^2N),
\end{eqnarray*}
which completes the proof. \qed

\section*{Acknowledgments}

The author would like to thank Maciej Radziejewski for helpful discussion and Pieter Moree for his remarks.


\begin{thebibliography}{gg}
\bibitem{BakerHarman} R.C. Baker and G. Harman, \emph{Shifted primes without large prime factors}, Acta Arith. \textbf{83} (1998), 331--361.
\bibitem{Carlitz-BinaryTerms} L. Carlitz, \emph{The number of terms in the cyclotomic polynomial $F_{pq}(x)$}, Amer. Math. Monthly \textbf{73} (1966), 979--981.
\bibitem{Hildebrand} A. Hildebrand, \emph{On a conjecture of Balog}, Proc. Amer. Math. Soc. \textbf{95} (1985), 517--523.
\bibitem{Hooley} C. Hooley, \emph{On the largest prime factor of $p+a$}, Mathematica \textbf{20} (1973), 135--143.
\bibitem{Lenstra-Vanishing} H.W. Lenstra,\emph{ Vanishing sums of roots of unity}, Proceedings, Bicentennial Congress Wiskundig Genootschap, Vrije Univ., Amsterdam (1978), Part II (1979), 249--268.
\bibitem{Migotti-BinaryBounds} A. Migotti, \emph{Aur Theorie der Kreisteilungsgleichung}, Z. B. der Math.-Naturwiss, Classe der Kaiserlichen Akademie der Wissenschaften, Wien, \textbf{87} (1883), 7--14.
\end{thebibliography}
\end{document}